\documentclass{amsart}
\usepackage{mathrsfs}
\usepackage{txfonts}
\usepackage{amssymb,amscd,amsthm}
\usepackage{latexsym}
\usepackage{hyperref}
\date{\today}

\newcommand{\Z}{{\mathbb Z}}
\newcommand{\R}{{\mathbb R}}

\newcommand{\e}{\varepsilon}

\newtheorem{theorem}{Theorem}[section]
\newtheorem{remark}[theorem]{Remark}
\newtheorem{lemma}[theorem]{Lemma}
\newtheorem{prop}[theorem]{Proposition}

\newtheorem{definition}[theorem]{Definition}
\newtheorem{claim}[theorem]{Claim}

\begin{document}
\title[Limit-Periodic Schr\"odinger Operators]{Limit-Periodic Schr\"{o}dinger Operators in the Regime of Positive Lyapunov Exponents }

\author[D.\ Damanik]{David Damanik}

\address{Department of Mathematics, Rice University, Houston, TX~77005, USA}

\email{damanik@rice.edu}

\urladdr{www.ruf.rice.edu/$\sim$dtd3}

\author[Z.\ Gan]{Zheng Gan}

\address{Department of Mathematics, Rice University, Houston, TX~77005, USA}

\email{zheng.gan@rice.edu}

\urladdr{math.rice.edu/$\sim$zg2}

\thanks{D.\ D.\ \& Z.\ G.\ were supported in part by NSF grant
DMS--0800100.}

\begin{abstract}
We investigate the spectral properties of discrete
one-dimensional Schr\"{o}dinger operators whose potentials are
generated by continuous sampling along the orbits of a minimal
translation of a Cantor group. We show that for given Cantor group
and minimal translation, there is a dense set of continuous sampling
functions such that the spectrum of the associated operators has
zero Hausdorff dimension and all spectral measures are purely
singular continuous. The associated Lyapunov exponent is a
continuous strictly positive function of the energy. It is possible
to include a coupling constant in the model and these results then
hold for every non-zero value of the coupling constant.
\end{abstract}

\maketitle

\section{Introduction}

This paper is a part of a sequence of papers devoted to the study of
spectral properties of discrete one-dimensional limit-periodic
Schr\"odinger operators. The first paper in this sequence,
\cite{dzg}, contains results in the regime of zero Lyapunov
exponents, while the present paper investigates the regime of
positive Lyapunnov exponents. Our general aim is to exhibit as rich
a spectral picture as possible within this class of operators. In
particular, we want to show that all basic spectral types are
possible and, in addition, in the case of singular continuous
spectrum and pure point spectrum, we are interested in examples with
positive Lyapunov exponents and examples with zero Lyapunov
exponents. From this point of view, the present paper will, to the
best of our knowledge for the first time, exhibit limit-periodic
Schr\"odinger operators with purely singular continuous spectrum and
positive Lyapunov exponents (whereas \cite{dzg} had the first
examples of limit-periodic Schr\"odinger operators with purely
singular continuous spectrum and zero Lyapunov exponents). Examples
with purely absolutely continuous spectrum have been known for a
long time, dating back to works of Avron and Simon \cite{as},
Chulaevsky \cite{c}, and Pastur and Tkachenko \cite{pt1,pt2} in the
1980's. These examples (must) have zero Lyapunov exponents. Examples
with pure point spectrum (and positive Lyapunov exponents at least
at many energies in the spectrum) can be found in P\"oschel's paper
\cite{p}; compare also the work of Chulaevsky and Molchanov
\cite{mc} (who have examples with zero Lyapunov exponents). In the third paper of this sequence we use P\"oschel's general theorem from \cite{p} to construct limit-periodic examples
with uniform pure point spectrum within our framework (actually these
examples have uniform localization of eigenfunctions); see \cite{dzg2}.

Our study is motivated by the recent paper \cite{a}, in which Avila
disproves a conjecture raised by Simon; see \cite[Conjecture
8.7]{s}. That is, he has shown that it is possible to have ergodic
potentials with uniformly positive Lyapunov exponents and
zero-measure spectrum. The examples constructed by Avila are
limit-periodic. In fact, the paper \cite{a} proposes a novel way of
looking at limit-periodic potentials. In hindsight, this way is
quite natural and provides one with powerful technical tools.
Consequently, we feel that a general study of limit-periodic
Schr\"odinger operators may be based on this new approach and we
have implemented this in \cite{dzg,dzg2} and the present paper. We anticipate that further
results may be obtained along these lines.

It has been understood since the early papers on limit-periodic
Schr\"odinger operators, and more generally almost periodic
Schr\"odinger operators, that these operators belong naturally to
the class of ergodic Schr\"odinger operators, where the potentials
are obtained dynamically, that is, by iterating an ergodic map and
sampling along the iterates with a real-valued function; see
\cite{cl,c2,pf} for general background. Indeed, taking the closure
in $\ell^\infty$ of the set of translates of an almost periodic
function on $\Z$ (i.e., the \emph{hull} of the function), one
obtains a compact Abelian group with a unique translation
invariant probability measure (Haar measure). In particular, the
shift on the hull is ergodic with respect to Haar measure and each
element of the hull may be obtained by continuous sampling (using
the evaluation at the origin, for example).

As pointed out by Avila, it is quite natural to take this one step
further. That is, once a compact Abelian group and a minimal
translation have been fixed, one is certainly not bound to sample
along the orbits merely with functions that evaluate a sequence at
one point. Rather, every continuous real-valued function on the
group is a reasonable sampling function. While this is quite
standard in the quasi-periodic case, we are not aware of any
systematic use of it in the context of limit-periodic potentials
before Avila's work \cite{a}.

The ability to fix the base dynamics and independently vary the
sampling functions is very useful in constructing examples of
potentials and operators that exhibit a certain desired spectral
feature. This has been nicely demonstrated in \cite{a} and is also
the guiding principle in our present work. As mentioned above, our
main motivation is to find examples of limit-periodic
Schr\"odinger operators with prescribed spectral type. From this
point of view, the singular continuity result we prove here is the
main result of the paper. However, there was additional motivation
to improve the zero measure result of Avila to a zero Hasudorff
dimension result. Recent work of Damanik and Gorodetski \cite{dg,dg2} focused on the weakly coupled Fibonacci Hamiltonian. This is an ergodic model that is not (uniformly) almost periodic.
Among the results obtained in \cite{dg,dg2}, there is a theorem that
states that the Hausdorff dimension of the spectrum, as a function
of the coupling constant, is continuous at zero. That is, as the
coupling constant approaches zero, the Hausdorff dimension of the
spectrum approaches $\dim_H([-2,2]) = 1$. When presenting this
result, the authors of \cite{dg,dg2} were asked whether this is a
universal feature, which holds for all potentials. Thus, our
purpose here is to show that there are indeed limit-periodic
potentials such that continuity at zero coupling fails in the
worst way possible, that is, the Hausdorff dimension of the
spectrum is zero for all non-zero values of the coupling
constant.\footnote{Our work was carried out right after the
preprint leading to the publication \cite{a} had been released.
That version proved zero-measure and did not discuss the Hausdorff
dimension issue. After we informed Avila about our results, we
learned from him that he had added a remark to the final version
of \cite{a} stating that a suitable modification of his proof of
zero measure yields zero Hausdorff dimension; see \cite[Remark
1.1]{a}.}

Let us now describe the models and results in detail. We consider discrete one-dimensional ergodic Schr\"{o}dinger operators
acting in $\ell^2(\Z)$ given by
\begin{equation}\label{oper}
[H^{\omega}_{f,T} \psi](n) = \psi(n+1) + \psi(n-1) + V_\omega(n) \psi(n)
\end{equation}
with
$$
V_\omega (n) = f(T^n (\omega)),
$$
where $\omega$ belongs to a compact space $\Omega$, $T: \Omega \rightarrow \Omega$ is a homeomorphism preserving an
ergodic Borel probability measure $\mu$ and $f : \Omega \to \R$ is a continuous sampling function. It is often beneficial
to study the operators $\{H^{\omega}_{f,T}\}_{\omega \in \Omega}$ as a family, as opposed to a collection of individual
operators, since the spectrum and the spectral type of $H^{\omega}_{f,T}$ are always $\mu$-almost surely independent of
$\omega$ due to ergodicity. Moreover, if $T$ is in addition minimal (i.e., all $T$-orbits are dense), then both the
spectrum and the absolutely continuous spectrum of $H^{\omega}_{f,T}$ are independent of $\omega$.

The Lyapunov exponent is defined as
\begin{align}\label{lya}
L(E,T,f) = \lim_{n\rightarrow \infty} \frac{1}{n} \int_{\Omega}
\log \| A^{(E,T,f)}_n(\omega) \| \, d\mu(\omega),
\end{align}
where $E \in \R$ is called the energy and $A^{(E,T,f)}_n$ is the $n$-step transfer matrix of \eqref{oper} defined as
\begin{equation}\label{equ:tran}
A^{(E,T,f)}_n(\omega) = S_{n-1} \dots S_0,\ \mathrm{where}\ S_i=\begin{pmatrix} E-f(T^i(\omega))& -1\\ 1& 0
\end{pmatrix}.
\end{equation}
By the Ishii-Pastur-Kotani theorem, the almost sure absolutely continuous spectrum of $H^{\omega}_{f,T}$ is given by the
essential closure of the set of energies where the Lyapunov exponent vanishes.

Next we make the spaces and homeomorphisms of especial interest to us explicit.

\begin{definition}
$\Omega$ is called a Cantor group if it is an infinite totally disconnected compact Abelian topological group.
\end{definition}

\begin{definition}
Let $\Omega$ be a Cantor group. For $\omega_1 \in \Omega$, let $T: \Omega \rightarrow \Omega$ be the translation by $\omega_1$, that is,
$T(\omega) = \omega_1 \cdot \omega$. $T$ is called minimal if $\{T^{n}(\omega) : n \in \Z\}$ is dense in $\Omega$ for every $\omega
\in \Omega$.
\end{definition}

We will restrict our attention to the case where $\Omega$ is a Cantor group and $T$ is a minimal translation. As mentioned
above, the operators $H^\omega_{f,T}$ have a common spectrum which we will denote by $\Sigma(f)$.

Here is our main result:

\begin{theorem} \label{t.main}
Suppose $\Omega$ is a Cantor group and $T$ is a minimal translation on $\Omega$. Then these exists a dense set $\mathcal{F} \subset
C(\Omega,\R)$ such that for every $f \in \mathcal{F}$ and every $\lambda \neq 0$, the following statements hold true:
$\Sigma(\lambda f)$ has zero Hausdorff dimension, $H^\omega_{\lambda f,T}$ has purely singular continuous spectrum for
every $\omega \in \Omega$, and $E \mapsto L(E,T,\lambda f)$ is a positive continuous function.
\end{theorem}

The proof of this theorem is based on the constructions in \cite{a}. We make several modifications to these constructions
to better control the size of the spectrum and to ensure that the potentials we construct are Gordon potentials. The latter
property then implies the absence of point spectrum, which in turn yields singular continuity since the absence of
absolutely continuous spectrum already follows from zero measure spectrum.

Let us state the Gordon property as a separate result.

\begin{definition}\label{d.gordon}
A bounded map $V: \Z \rightarrow \R$ is called a Gordon potential if there exist positive integers $q_i \rightarrow \infty$
such that
$$
\max_{1\leq n \leq q_i}\left|V(n)-V(n\pm q_i) \right|\leq i^{-q_i}
$$
for every $i \ge 1$.
\end{definition}

Clearly, if $V$ is a Gordon potential, so is $\lambda V$ for every $\lambda \in \R$. A key part in proving Theorem~\ref{t.main}
is to establish the following result:

\begin{theorem} \label{t.gordon}
Suppose $\Omega$ is a Cantor group. Then these exists a dense set $\mathcal{F} \subset C(\Omega,\R)$ such that for every $f \in
\mathcal{F}$, every minimal translation $T : \Omega \rightarrow \Omega$, every $\omega \in \Omega$, and every $\lambda \neq 0$,
$\lambda f(T^n(\omega))$ is a Gordon potential.
\end{theorem}

\section{Preliminaries}

\subsection{Hausdorff Measures and Dimensions}

For our relatively restricted purposes, we will simply recall the definition of Hausdorff measures and Hausdorff dimension in this subsection.
We refer the reader to \cite{r} for more information.

\begin{definition}
Let $A \subseteq \R$ be a subset. A countable collection of intervals $\{ b_n \}^{\infty}_{n=1}$ is called a $\delta$-cover of $A$ if $A \subset \bigcup^{\infty}_{n=1} b_n$ with $|b_n| < \delta$ for all $n$'s. {\rm (}Here, $|\cdot|$ denotes Lebesgue measure, and we will
adopt this notation throughout the paper.{\rm )}
\end{definition}

\begin{definition}
Let $\alpha \in \R$. For any subset $A \subseteq  \R$, the $\alpha$-dimensional Hausdorff measure of $A$ is defined as
\begin{equation}\label{equ:haus}
h^{\alpha}(A) = \lim_{\delta \rightarrow 0} \inf_{\text{$\delta$-covers}}
\sum^{\infty}_{n=1} |b_n|^{\alpha} .
\end{equation}
\end{definition}

The quantity $h^\alpha(A)$ is well defined as an element of $[0,\infty]$ since $\inf_{\text{$\delta$-covers}}
\sum^{\infty}_{n=1} |b_n|^{\alpha} $ is monotonically increasing as $\delta$ decreases to zero and therefore
the limit in \eqref{equ:haus} exists. Restricted to the Borel sets, $h^1$ coincides with Lebesgue measure and
$h^0$ is the counting measure. If $\alpha < 0$, we always have $h^{\alpha}(A) = \infty$ for any $A \neq \emptyset$, while if $\alpha > 1$, $h^{\alpha}(\R) = 0$.

It is not hard to see that for every $A \subseteq \R$, there is a unique $\alpha \in [0,1]$, called the Hausdorff
dimension $\dim_H(A)$ of $A$, such that $h^{\beta}(A) = \infty$ for every $\beta < \alpha$ and $h^{\beta}(A) = 0$
for every $\beta > \alpha$. In particular, every $A \subseteq \R$ with $|A| > 0$ must have $\dim_H(A) = 1$.

\subsection{Minimal Translations of Cantor Groups and Limit-Periodic Potentials}

In this subsection we recall how the one-to-one correspondence between hulls of limit-periodic sequences and potential families
generated by minimal translations of Cantor groups and continuous sampling functions exhibited by Avila in \cite{a} arises.

\begin{definition}
Let $S : \ell^\infty(\Z) \to \ell^\infty(\Z)$ be the shift operator, $(S V)(n) = V(n+1)$. A two-sided sequence $V \in
\ell^\infty(\Z)$ is called periodic if its $S$-orbit is finite and it is called limit-periodic if it belongs to the closure of the
set of periodic sequences. If $V$ is limit-periodic, the closure of its $S$-orbit is called the hull and denoted by
$\mathrm{hull}_V$.
\end{definition}

The first lemma (see \cite[Lemma~2.1]{a}) shows how one can write the elements of the hull of a limit-periodic function in the form
\begin{equation}\label{pot}
V_\omega (n) = f(T^n (\omega)) , \quad \omega \in \Omega, \; n \in
\Z
\end{equation}
with a minimal translation $T$ of a Cantor group and a sampling function $f \in C(\Omega,\R)$:

\begin{lemma}
Suppose $V$ is limit-periodic. Then, $\Omega : = \mathrm{hull}_V$ is compact and has a unique topological group structure with
identity $V$ such that $\Z \ni k \mapsto S^k V \in \mathrm{hull}_V$ is a homomorphism. Moreover, the group structure
is Abelian and there exist arbitrarily small compact open neighborhoods of $V$ in $\mathrm{hull}_V$ that are finite index
subgroups.
\end{lemma}

In particular, $\Omega = \mathrm{hull}_V$ is a Cantor group, $T = S|_\Omega$ is a minimal translation, and every element of $\Omega$
may be written in the form \eqref{pot} with the continuous function $f(\omega) = \omega(0)$.

The second lemma (see \cite[Lemma~2.2]{a}) addresses the converse:

\begin{lemma} \label{lem:avila}
Suppose $\Omega$ is a Cantor group, $T : \Omega \to \Omega$ is a
minimal translation, and $f \in C(\Omega,\R)$. Then, for every
$\omega \in \Omega$, the element $V_\omega$ of $\ell^\infty(\Z)$
defined by \eqref{pot} is limit-periodic and we have
$\mathrm{hull}_{V_\omega} = \{ V_{\tilde \omega} \}_{\tilde \omega
\in \Omega}$.
\end{lemma}

These two lemmas show that a study of limit-periodic potentials can be carried out by considering potentials of the form
\eqref{pot} with a minimal translation $T$ of a Cantor group $\Omega$ and a continuous sampling function $f$. As shown for the
first time in the context of limit-periodic potentials by Avila in \cite{a}, it is often advantageous to fix $\Omega$ and $T$ and to
vary $f$.

\subsection{Periodic Sampling Functions, Potentials, and Schr\"odinger Operators}

In this subsection we discuss the periodic case. For example, which sampling functions $f \in C(\Omega,\R)$ give rise
to periodic potentials for some or all $(\omega,T)$? Moreover, what can then be said about the associated Schr\"odinger operators?

\begin{definition}
Suppose $\Omega$ is a Cantor group and $T : \Omega \to \Omega$ is a minimal translation. We say that a sampling function
$f \in C(\Omega,\R)$ is $n$-periodic with respect to $T$ if $f(T^n(\omega)) = f(\omega)$ for every $\omega \in \Omega$.
\end{definition}

\begin{prop}\label{inde}
Let $f \in C(\Omega,\R)$. If $f(T^{n_0+m}(\omega_0)) =
f(T^{m}(\omega_0))$ for some $\omega_0 \in \Omega$, some minimal
translation $T : \Omega \to \Omega$ and every $m \in \Z$, then for
every minimal translation $\tilde{T} : \Omega \to \Omega$, $f$ is
$n_0$-periodic with respect to $\tilde{T}$.
\end{prop}

\begin{proof}
Let $\varphi : \Omega \rightarrow \ell^{\infty}(\Z)$,
$\varphi(\omega) = (f(T^{n}(\omega)))_{n \in \Z}$. Since $T$ is
minimal, the closure of $\{T^{n}(\omega_0) : n \in \Z\}$ is
$\Omega$. By Lemma~\ref{lem:avila} we have $ \varphi(\Omega) =
\mathrm{hull}(\varphi(\omega_0)).$ Since $f(T^{n_0 + m}(\omega_0)) =
f(T^{m}(\omega_0))$ for any $m \in \Z$,
$\mathrm{hull}(\varphi(\omega_0))$ is a finite set. Then for any
$\omega \in \Omega$, $(f(T^{n}(\omega)))_{n \in \Z}$ is some element
in $\mathrm{hull}(\varphi(\omega_0))$. Since every element in
$\mathrm{hull}(\varphi(\omega_0))$ is $n_0$-periodic,
$(f(T^{n}(\omega)))_{n \in \Z}$ is $n_0$-periodic. This shows that
$f$ is $n_0$-periodic with respect to $T$. That is, we have
$f(T^{n_0+m}(\omega)) = f(T^{m}(\omega))$ for every $\omega \in
\Omega$ and $m \in \Z$.

Assume $T$ is the minimal translation by $\omega_1$ and let
$\tilde{T}$ be another minimal translation by $\omega_2$. By the
previous analysis, we have $f(\omega^{n_0 + m}_1\cdot\omega) =
f(\omega^{m}_1\cdot\omega)$ for every $m \in \Z$ and every $\omega
\in \Omega$. If $\omega_2$ is equal to $\omega^{q}_1$ for some
integer $q$, obviously we have $f(\tilde{T}^{n_0}(\omega)) =
f((\omega^{q}_1)^ {n_0}\cdot\omega) = f(\omega)$ for any $\omega
\in \Omega$. If not, since $\{\omega^{n}_1 : n \in \Z \}$ is dense
in $\Omega$ (this follows from the minimality of $T$), we have
$\lim_{k\rightarrow \infty} \omega^{n_k}_1 = \omega_2$, and then
$f(\omega^{n_0}_2\cdot\omega) = \lim_{k\rightarrow \infty}
f((\omega^{n_k}_1)^{n_0}\cdot \omega) = f(\omega)$. The result
follows.
\end{proof}

The above proposition tells us that the periodicity of $f$ is
independent of $T$. Thus we may say $f$ is $n$-periodic without
making a minimal translation explicit.

Next we recall from \cite{a} how periodic sampling functions in
$C(\Omega,\R)$ can be constructed. Given a Cantor group $\Omega$,
a compact subgroup $\Omega_0$ with finite index (such subgroups
can be found in any neighborhood of the identity element; see
above), and $f \in C(\Omega,\R)$, we can define a periodic
$f_{\Omega_0} \in C(\Omega,\R)$ by
$$
f_{\Omega_0}(\omega) = \int_{\Omega_0} f(\omega \cdot \tilde \omega) \, d\mu_{\Omega_0}(\tilde \omega).
$$
Here, $\mu_{\Omega_0}$ denotes Haar measure on $\Omega_0$. This
shows that the set of periodic sampling functions is dense in
$C(\Omega,\R)$. Moreover, as already noted in \cite{a}, there
exists a decreasing sequence of Cantor subgroups $\Omega_k$ with
finite index $n_k$ such that $\bigcap \Omega_k = \{e\}$, where $e$
is the identity element of $\Omega$. Let $P_k$ be the set of
sampling functions defined on $\Omega / \Omega_k$, that is, the
elements in $P_k$ are $n_k$-periodic potentials. Denote by $P$ the set of all periodic sampling functions. Then, we have $P_{k} \subset P_{k+1}$ (which implies $n_k \mid n_{k+1}$) and $P = \bigcup P_k$.

\begin{prop}\label{indely}
Let $f$ be $p$-periodic. Then, for every $\omega \in \Omega$,
\begin{align} \label{lya2}
L(E,T,f) & = \lim_{m \rightarrow \infty} \frac{1}{m}  \log \|
A^{(E,T,f)}_m (\omega) \| \\
\nonumber & = \frac{1}{p} \log \rho(A^{(E,T,f)}_{p} (e)),
\end{align}
where $\rho(A^{(E,T,f)}_{p} (e))$ is the spectral radius of $A^{(E,T,f)}_{p} (e)$. In particular, if
restricted to periodic sampling functions, the Lyapunov exponent is
a continuous function of both the energy $E$ and the sampling
function.
\end{prop}

\begin{proof}
If $f$ is $p$-periodic, as in the proof of Proposition \ref{inde},
for every $\omega$, $(f({T}^{n}(\omega)))_{n \in \Z}$ is some
element of the orbit of $(f(T^{n}(e)))_{n \in \Z}$, and so its
monodromy matrix (i.e., the transfer matrix over one period) is a
cyclic permutation of the monodromy matrix associated with
$f({T}^{n}(e))$. Thus $\mathrm{Tr} A^{(E,T,f)}_{p}(\omega)$ is
independent of $\omega$, and since $\det A^{(E,T,f)}_{p}(\omega) = 1
$, we can conclude that the eigenvalues of $A^{(E,T,f)}_{p}(\omega)$
are independent of $\omega$. So the logarithm of the spectral radius
of $A^{(E,T,f)}_{p}(\omega)$ is independent of $\omega$ and
\eqref{lya2} follows. The continuity statement follows readily.
\end{proof}

\begin{lemma}\label{lem:limitlya}
Let $f_n \in C(\Omega,\R)$ be a sequence of periodic sampling functions converging uniformly to $f_\infty \in C(\Omega,\R)$. Assume $
\lim_{n\rightarrow \infty} L(E,T,f_n)$ exists for every $E$ and the convergence is uniform. Then we have that $L(E,T,f_\infty)$
coincides with $\lim_{n\rightarrow \infty} L(E,T,f_n)$ everywhere.
\end{lemma}

\begin{proof}
Since $\lim_{n\rightarrow \infty} L(E,T,f_n)$ exists everywhere,
from \cite[Lemma 2.5]{a}, we have $L(E,T,f_n) \rightarrow
L(E,T,f_\infty)$ in $L^{1}_{loc}.$ So $L(E,T,f_\infty)$ coincides
with $\lim_{n\rightarrow \infty} L(E,T,f_n)$ almost everywhere.
From Proposition \ref{indely}, $L(E,T,f_n)$ is a continuous
function, and by uniform convergence, we have that
$\lim_{n\rightarrow \infty} L(E,T,f_n)$ is also a continuous
function. Since $L(E,T,f_\infty)$ is a subharmonic function
(cf.~\cite[Theorem 2.1]{cs}), we get that $L(E,T,f_\infty) =
\lim_{n\rightarrow \infty} L(E,T,f_n)$ for every $E$. The
statement follows.
\end{proof}

To conclude this subsection on the periodic case, we state two
lemmas. The first is well known and the second is \cite[Lemma
2.4]{a}.

\begin{lemma} \label{lem:ab}
Let $f \in C(\Omega,\R)$ be $p$-periodic. \\
(i). The spectrum of $H^\omega_{f,t}$ is purely absolutely continuous for every $\omega \in \Omega$ and $\Sigma(f)$ is made of
$p$ bands {\rm (}compact intervals whose interiors are disjoint{\rm )}. \\
(ii). $\Sigma(f) = \{E\in \R : L(E,T,f) = 0\}.$
\end{lemma}

\begin{lemma}\label{lem:measure}
Let $f \in C(\Omega,\R)$ be $p$-periodic. \\
(i). The Lebesgue measure of each band of $\Sigma(f)$ is at most $\frac{2\pi}{p}$.\\
(ii). Let $C\ge 1$ be such that for every $E \in \Sigma(f)$, there
exist $\omega\in \Omega$ and $k\ge 1$ such that $\|
A_k^{(E,T,f)}(\omega) \| \ge C$. Then, $| \Sigma(f) | \le \frac{4\pi p}{C}$.
\end{lemma}

\section{Proof of the Theorems}

Assume $\Omega$ and $T$ are given. For convenience, we write
$A^{(E,f)}_n(\omega) = A^{(E,f,T)}_n(\omega)$, $A^{(E,f)}_n =
A^{(E,f,T)}_n(e)$, and $L(E,f) = L(E,T,f)$. Since $T:\Omega
\rightarrow \Omega$ is a minimal translation, the homomorphism $\Z
\rightarrow \Omega$, $n \rightarrow T^ne$ is injective with dense
image in $\Omega$, and we can write $f(n)=f(T^n(e))$ without any
conflicts.

We need two more lemmas before proving our theorems. More
precisely, we will make further use of the constructions which
play central roles in the proof of these two lemmas.

\begin{lemma}\label{lem:a1}
Let $B$ be an open ball in $C(\Omega,\R)$, let $F \subset P \cap B$ be finite, and let $0 < \e < 1$. Then there
exists a sequence $F_K \subset P \cap B$ such that \\ [1mm] (i). $L(E ,\lambda F_K) > 0$ whenever $\e \leq |\lambda| \leq \e^{-1}$,
$E \in \R$,\\
(ii). $L(E, \lambda F_K) \rightarrow L(E,\lambda F)$ uniformly on compacts (as functions of $(E,\lambda) \in \R^2$).
\end{lemma}

This is \cite[Lemma 3.1]{a}. As in \cite{a}, we use the notation
$$
L(E, \lambda F) = \frac{1}{\# F} \sum_{f \in F} L(E,T, \lambda f),
$$
where $F$ is a finite family of sampling functions (with
multiplicities!) and $\lambda \in \R$. The proof of this lemma is
constructive. We will describe this construction explicitly in the
proof of Theorem~\ref{t.main} for the reader's convenience.

\begin{lemma}\label{lem:a2}
Suppose $B$ is an open ball in $C(\Omega,\R)$ and $F \subset P \cap B$ is a finite family of sampling functions. Then for every
$N \ge 2$ and $K$ sufficiently large, there exists $F_K \subset P_K \cap B $ such that \\[1mm]
(i). $L(E, \lambda F_K) \rightarrow L(E, \lambda F)$ uniformly on compacts (as functions of $(E,\lambda) \in \R^2$).\\
(ii). The diameter of $F_K$ is at most $n_K^{-N/2}$.
\end{lemma}

This lemma is a variation of \cite[Lemma 3.2 ]{a}. We will prove
this lemma using suitable modifications of Avila's arguments. Some
of these modifications, which will later enable us to prove the
Gordon property, are not apparent from the statement of the lemma.
We will give detailed arguments for the modified parts of the
proof and refer the reader to \cite{a} for the parts that are
analogous.

\begin{proof}[Proof of Lemma \ref{lem:a2}] Assume that $F = \{f_1,f_2,\dots,f_m\}
\subset C(\Omega,\R)$ is a finite family of $n_k$-periodic sampling
functions with $n_k \ge 2$, and let $K > k$ be large enough. We
construct $F^{\vec{t}}_K$ as follows. Let $n_K = m n_k r + d,\ 0
\leq d \leq m n_k - 1$. Let $I_j = [jn_k,(j+1)n_k -1] \subset \Z$
and let $0 = j_0 < j_1 < \dots < j_{m-1} < j_{m} = n_K/n_k$ be a
sequence such that $j_{i+1} - j_i = r + 1$ when $0 \leq i < d/n_k$ and $j_{i+1} - j_{i} = r$ when $d/n_k
\leq i \leq m - 1$. Define an $n_K$-periodic $f$ as follows. For $0 \leq l \leq n_K - 1$, let $j$ be such that $l \in I_j$ and let $i$ be such that $j_{i-1} \leq j < j_i$ and then let $f(l) = f_i(l)$.
Next, for any sequence $\vec{t} = (t_1,t_2,\dots,t_{m})$ with $t_i
\in \{0,1,\dots,r - 1\}$, we define an $n_K$-periodic $f^{\vec{t}}_K
$ as follows. If $j = j_i - 1$ for some $1 \leq i < m$, we let
$f^{\vec{t}}_K(l) = f(l) + r^{-N} t_i$, and if $j = j_m -2$, we let
$f^{\vec{t}}_K(l) = f(l) + r^{-N} t_m$. Otherwise we let
$f^{\vec{t}}_K(l) = f(l)$. Let $F^{\vec{t}}_K$ be the family
consisting of all $f^{\vec{t}}_K$'s. The statement (ii) is clear for
large $K$. (Note: in \cite{a}, Avila's construction is such that if
$j = j_i - 1$ for some $1 \leq i \leq m$, then $f^{\vec{t}}_K(l) =
f(l) + r^{-20} t_i$; otherwise, $f^{\vec{t}}_K(l) = f(l)$.)

For fixed $E$ and $\lambda$, we let $A^{(E,\lambda
f^{\vec{t}}_K)}_{n_K} = C^{(t_m,m)}B^{(m)} \dots
C^{(t_1,1)}B^{(1)}$, where $B^{(i)} = (A^{(E,\lambda
f_i)}_{n_k})^{j_i - j_{i-1} - 1}, 1 \leq i \leq m-1$ and $B^{(m)}
= (A^{(E,\lambda f_m)}_{n_k})^{j_m - j_{m-1} - 2}$, and
$C^{(t_i,i)} = A^{(E - \lambda r^{-N}t_i,\lambda f_i)}_{n_k}, 1
\leq i \leq m-1$ and $C^{(t_m,m)} = A^{(E,\lambda f_m)}_{n_k}
A^{(E - \lambda r^{-N}t_m,\lambda f_m)}_{n_k}$. When $E$ and
$\lambda$ are in a compact set, the norm of the $C^{(t_i,i)}$-type
matrices is bounded as $r$ grows, while the norm of the
$B^{(i)}$-type matrices may get large.

Notice that our perturbation here is $r^{-N}t$ (as opposed to Avila's $r^{-20}t$
perturbation in \cite[Lemma 3.2 ]{a}), so \cite[Claim 3.7]{a} should be replaced by the following version:\\[1mm]
\indent``Let $s_j$ be the most contracted direction of $\hat{B}^{(j)}$ and let $u_j$
be the image under $\hat{B}^{(j)}$ of the most expanded direction. Call $\vec{t}$ $j$-nice, $1 \leq j \leq d$,
if the angle between $\hat{C}^{(j)}u_j$ and $S_{j+1}$ (less than $\pi$) is at least $r^{-3 N}$ with the
convention that $j + 1 = 1$ for $j = d$. Let $r$ be sufficiently large, and let $\vec{t}$ be $j$-nice.
If $z$ is a non-zero vector making an angle at least $r^{-4 N}$ with $s_j$, then $z^{'} = \hat{C}^{(j)}
\hat{B}^{(j)} z$ makes an angle at least $r^{-4 N}$ with $S_{j+1}$ and $\| z^{'} \| \ge \| \hat{B}^{(j)} \| r^{-5 N} \| z \|.$"\\[1mm]
The proof of \cite[Claim 3.7]{a} can be applied to get the above
version of the claim with the corresponding quantitative
modification. Moreover, we have also made a little shift in the
perturbation, so $C^{(t_m,m)} = A^{(E,\lambda f_m)}_{n_k} A^{(E -
\lambda r^{-N}t_m,\lambda f_m)}_{n_k}$, while Avila's $C^{(t_m,m)}
= A^{(E - \lambda r^{-20}t_m,\lambda f_m)}_{n_k}$. \cite[Claim
3.8]{a} still holds, but Avila's proof of \cite[Claim 3.8]{a}
cannot be applied directly. To this end we prove the following
claim:

\begin{claim} \label{matrix}
For any $M \in \mathrm{SL}(2,R)$, there are $m_1, m_2 \in
(0,\infty)$ with the following property. Suppose $A$ and $B$ are two
vectors in $\R^{2}$, and $\Delta \theta$ is the angle between $A$
and $B$ with $0 < \Delta \theta \leq \pi$. Let $\Delta
\tilde{\theta}$ be the angle between $M A$ and $M B$ {\rm (}again so
that $0 < \Delta \tilde{\theta} \leq \pi${\rm )}. Then, $ m_1 \Delta
\theta \leq \Delta \tilde{\theta} \leq m_2 \Delta \theta. $
\end{claim}

\begin{proof}
By the singular value decomposition (see \cite[Theorem 2.5.1]{svd}),
there exist $O_1$ and $O_2$ in $\mathrm{SO}(2,R)$ such that $ M =
O_1  S  O_2,$ where $S$ is a diagonal matrix. Since $O_1$ and $O_2$
are rotations on $\R^2$, it is sufficient to consider
$$
S = \begin{pmatrix} \mu_1 & 0\\
0&\mu^{-1}_1
\end{pmatrix}.
$$

Without loss of generality, assume $\mu_1 \ge 1$. Let $A =
(a,b)^{t}$ ($t$ denotes the transpose of vectors) and $B =
(c,d)^{t}$ be two normalized vectors, and let $\theta_{A}$ and
$\theta_{B}$ be the argument of $A$ and the argument of $B$
respectively. Let $\tilde{A} = S A = ( a\mu_1, b/\mu_1 )^t$ with
the argument $\theta_{\tilde{A}}$ and $\tilde{B} = S B = ( c\mu_1,
d/\mu_1 )^t$ with the argument $\theta_{\tilde{B}}$.

We adopt the following notation for convenience. Let
$\uppercase\expandafter{\romannumeral1}$,
$\uppercase\expandafter{\romannumeral2}$,$\uppercase\expandafter{\romannumeral3}$,$\uppercase\expandafter{\romannumeral4}$
denote one of two vectors in the first quadrant (including
$\{(x,0): x \ge 0 \}$), the second quadrant (including $\{(0,y): y
> 0 \}$), the third quadrant (including $\{(x,0): x < 0 \}$) and
the fourth quadrant (including $\{(0,y): y < 0 \}$), respectively.
Then
$(\uppercase\expandafter{\romannumeral1},\uppercase\expandafter{\romannumeral1})$
denotes that both two vectors are in the first quadrant,
$(\uppercase\expandafter{\romannumeral1},\uppercase\expandafter{\romannumeral2})$
denotes that one vector is in the first quadrant while the other
is in the second quadrant, and so on.

We will need the following observation:
\begin{equation}\label{e.observe}
0 < \theta_1, \theta_2 < \pi/2 \text{ and } \tan \theta_1 \ge
\frac{\tan \theta_2}{\mu^2_1} \quad \Rightarrow \quad \theta_1 \ge
\frac{\theta_2}{4 \mu^2_1}.
\end{equation}
Indeed, since $\tan \theta_1 \ge \frac{1}{\mu^2_1} \tan \theta_2
\ge \frac{1}{2 \mu^2_1} \theta_2$ and $0 < \frac{\theta_2}{2
\mu^2_1} < 1$, we have
$$
\theta_1 \ge \arctan{ \frac{\theta_2}{2 \mu^2_1}} =
\frac{\theta_2}{2 \mu^2_1} - {(\frac{\theta_2}{2 \mu^2_1})^3}/3 +
O((\frac{\theta_2}{2 \mu^2_1})^5) \ge \frac{\theta_2}{4 \mu^2_1}.
$$

For the proof of Claim~\ref{matrix}, we consider two cases.\\[1mm]
\textit{Case 1}. $\pi/2 \leq \Delta\theta \leq \pi$.  Here $A$ and
$B$ cannot be in the same quadrant. Notice that the impact of $S$
on vectors is to move them closer to the $x$-axis and keep them in
the same quadrant. Thus, for the subcases
$(\uppercase\expandafter{\romannumeral1},\uppercase\expandafter{\romannumeral2})$,
$(\uppercase\expandafter{\romannumeral1},\uppercase\expandafter{\romannumeral3})$,
$(\uppercase\expandafter{\romannumeral2},\uppercase\expandafter{\romannumeral4})$
and
$(\uppercase\expandafter{\romannumeral3},\uppercase\expandafter{\romannumeral4})$,
we can easily conclude that $\Delta\theta/2 \leq
\Delta\tilde{\theta} \leq 2 \Delta\theta$. There are two subcases
left,
$(\uppercase\expandafter{\romannumeral1},\uppercase\expandafter{\romannumeral4})$
and
$(\uppercase\expandafter{\romannumeral2},\uppercase\expandafter{\romannumeral3})$.
We will discuss
$(\uppercase\expandafter{\romannumeral1},\uppercase\expandafter{\romannumeral4})$;
the method can be readily adapted to
$(\uppercase\expandafter{\romannumeral2},\uppercase\expandafter{\romannumeral3})$.
For
$(\uppercase\expandafter{\romannumeral1},\uppercase\expandafter{\romannumeral4})$,
if $\theta_{A} = 0$ and $\theta_{B} = 3\pi/2$, then
$\theta_{\tilde{A}}$ and $\theta_{\tilde{B}}$ are also $0$ and
$3\pi/2$ respectively, and so $\Delta \tilde{\theta} = \Delta
\theta$; if not, without loss of generality, assume that $A$ is in
the first quadrant with $\pi/4 \leq \theta_{A} < \pi/2 $, then $
\tan{\theta_{\tilde{A}}} = \frac{b}{a\mu^2_1} =
\frac{\tan{\theta_{A}}}{\mu^2_1}, $ and by \eqref{e.observe}, we
have
$$\Delta \tilde{\theta} \ge \theta_{\tilde{A}} \ge \frac{\theta_{A}}{4\mu^2_1}  \ge \frac{\Delta \theta}{16 \mu^2_1} $$
($ \theta_{A} \ge \Delta \theta /4 $ since $\theta_{A} \ge \pi /4$)
and then $\frac{\Delta\theta}{16 \mu^2_1} \leq \Delta\tilde{\theta}
\leq 2
\Delta\theta.$\\[1mm]
\textit{Case 2}. $0 < \Delta\theta < \pi/2$. In this case,
$(\uppercase\expandafter{\romannumeral1},\uppercase\expandafter{\romannumeral3})$
and
$(\uppercase\expandafter{\romannumeral2},\uppercase\expandafter{\romannumeral4})$
are impossible. We will divide the following proof into three parts.\\
(1). We discuss
$(\uppercase\expandafter{\romannumeral1},\uppercase\expandafter{\romannumeral1})$
here; the argument may be readily adapted to
$(\uppercase\expandafter{\romannumeral2},\uppercase\expandafter{\romannumeral2})$,
$(\uppercase\expandafter{\romannumeral3},\uppercase\expandafter{\romannumeral3})$,
and
$(\uppercase\expandafter{\romannumeral4},\uppercase\expandafter{\romannumeral4})$.
Without loss of generality, assume $\Delta\theta = \theta_{A} -
\theta_{B}$, then we get
$$
\tan{\Delta\tilde{\theta}} = \frac{\mu^2_1(bc - ad)}{bd + \mu^4_1
ac} \ge \frac{\tan{\Delta\theta}}{\mu^2_1},
$$
and by \eqref{e.observe}, we get $ \frac{\Delta\theta}{4 \mu^2_1}
\leq \Delta\tilde{\theta}$. Similarly, we will get
$\Delta\tilde{\theta} \leq {4 \mu^2_1} \Delta\theta$ since
$\tan{\Delta\tilde{\theta}} \leq \mu^2_1 \tan{\Delta\theta}$, and
so $ \frac{\Delta\theta}{4 \mu^2_1} \leq \Delta\tilde{\theta} \leq
{4 \mu^2_1} \Delta\theta$ follows.\\
(2). We discuss
$(\uppercase\expandafter{\romannumeral1},\uppercase\expandafter{\romannumeral4})$
here; an adaptation handles
$(\uppercase\expandafter{\romannumeral2},\uppercase\expandafter{\romannumeral3})$.
Without loss of generality, assume $\theta_A \ge \Delta\theta/2$.
Obviously, we have $\Delta\tilde{\theta} \leq \Delta{\theta}$.
Conversely, we have $\frac{\Delta{\theta}}{16 \mu^2_1} \leq
\Delta\tilde{\theta}$ (it is essentially the same as
$(\uppercase\expandafter{\romannumeral1},\uppercase\expandafter{\romannumeral4})$
in \textit{Case 1}), and so $\frac{\Delta{\theta}}{16 \mu^2_1}
\leq
\Delta\tilde{\theta} \leq \Delta{\theta}$ follows.\\
(3). We discuss
$(\uppercase\expandafter{\romannumeral1},\uppercase\expandafter{\romannumeral2})$
here, and the method can be applied to
$(\uppercase\expandafter{\romannumeral3},\uppercase\expandafter{\romannumeral4})$.
Obviously, we have $\Delta\theta \leq \Delta\tilde{\theta}.$
Without loss of generality, assume that $A$ is in the first
quadrant and makes an angle $h_A$ with the $y$-axis and that $B$
is in the second quadrant and makes an angle $h_B$ with the
$y$-axis. Clearly, $\Delta \theta = h_A + h_B$. Let
$h_{\tilde{A}}$ and $h_{\tilde{B}}$ be the angle between the
$y$-axis and $\tilde{A}$ and the angle between the $y$-axis and
$\tilde{B}$, respectively. By \eqref{e.observe}, we conclude that
$h_{\tilde{A}} \leq 4\mu^2_1 h_A$ since $\tan{h_{\tilde{A}}} =
\mu^2_1 \tan{h_A}$. Similarly, we get $h_{\tilde{B}} \leq 4\mu^2_1
h_B$. So it follows that $\Delta\theta \leq \Delta{\tilde{\theta}}
= h_{\tilde{A}} + h_{\tilde{B}} \leq 4\mu^2_1 (h_A + h_B) =
4\mu^2_1 \Delta
\theta.$\\[1mm]
\indent Through the above analysis, we see that
$\frac{\Delta\theta}{16 \mu^2_1} \leq \Delta\tilde{\theta} \leq
{16 \mu^2_1} \Delta\theta$, concluding the proof of
Claim~\ref{matrix}.
\end{proof}

By this claim, we can modify the last paragraph of the proof of \cite[Claim 3.8]{a} as stated below and then our lemma follows.
%(Note: we leave the necessary modifications to the corresponding parts of \cite[Proof of Lemma 3.2]{a} to the reader.)
\\[1mm]
\indent ``If $r$ sufficiently large, we conclude that for every $0 \leq l \leq r-2$, there exists a rotation $R_{l,j}$
of angle $\theta_j$ with $r^{-2.5N} < \theta_j < r^{-0.3N}$ such that $C^{(l+1,i_j)} u_j = R_{l,j} C^{(l,i_j)} u_j$.
It immediately follows that there exists at most one choice of $0 \leq t_{i_j} \leq r-1$ such that $C^{(t_{i_j},i_j)}u_j$
has angle at most $r^{-3N}$ with $s_{j+1}$, as desired."\\[1mm]
\indent We would like to explain how to obtain the statement described in the paragraph above. If $r$ is sufficiently
large, it is not hard to conclude that for every $0 \leq l \leq r-2$, there exists a rotation $\tilde {R}_{l,j}$ of angle
$\tilde{\theta}_j$ with $r^{-2N} < \tilde{\theta}_j < r^{-0.5N}$ such that $A^{(E - \lambda r^{-N}(l+1),\lambda f_{i_j})}_{n_k} u_j
= \tilde {R}_{l,j} A^{(E - \lambda r^{-N}l,\lambda f_{i_j})}_{n_k} u_j$. If $i_d = m$, we have
\begin{align}
 \label{angle}C^{(l+1,m)} u_m &= A^{(E,\lambda f_m)}_{n_k} A^{(E - \lambda r^{-N} (l+1),\lambda f_m)}_{n_k} u_m \\
& \nonumber = A^{(E,\lambda f_m)}_{n_k} \tilde{R}_{l,m} A^{(E -
\lambda r^{-N} (l),\lambda f_m)}_{n_k} u_m.
\end{align}
Since $A^{(E,\lambda f_m)}_{n_k} \in \mathrm{SL}(2,R)$ is independent of $r$, we can apply Claim~\ref{matrix}  to \eqref{angle} so that we
have
\begin{align*}
C^{(l+1,m)} u_m  &= A^{(E,\lambda f_m)}_{n_k} \tilde{R}_{l,m} A^{(E - \lambda r^{-N} (l),\lambda f_m)}_{n_k} u_m \\&= R_{l,m}
A^{(E,\lambda f_m)}_{n_k} A^{(E - \lambda r^{-N} (l),\lambda f_m)}_{n_k} u_m \\&= R_{l,m} C^{(l,m)} u_m,
\end{align*}
where $R_{l,m}$ is a rotation of angle $\theta_m$ with $r^{-2.5N} < \theta_m < r^{-0.3N}$. Then the above paragraph follows.
\end{proof}

Recall the definition of a Gordon potential given in Definition~\ref{d.gordon}. The importance of Gordon potentials lies in the following lemma,
which (in a slightly weaker form) first appeared in \cite{g}.

\begin{lemma}[{Gordon Lemma}]\label{lem:gordon}
Suppose $V$ is a Gordon potential. Then the Schr\"{o}dinger operator with potential $V$ has no eigenvalues.
\end{lemma}

Now we can give the

\begin{proof}[Proof of Theorem~\ref{t.main}]
Given a $p_0$-periodic $f \in C(\Omega,\R)$ and $0 < \e_0 < 1$,
consider $B_{\e_0}(f) \subset C(\Omega,\R)$. (We will work within
this ball. The denseness of periodic potentials then implies the
denseness of our constructed limit-periodic potentials.) Let $N$
from Lemma~\ref{lem:a2} be $2$. Let $\e_1 = \frac{\e_0}{10}$. By
Lemma~\ref{lem:a1}, there exists a finite family $F_1 = \{
f_1,f_2,\dots f_{m_1} \}$ of $p_1$-periodic sampling functions
such that $F_1 \subset B_{\e_0}(f)$ and $L(E,\lambda F_1) >
\delta_1$ for some $0 < \delta_1 < 1$ whenever $\e_1 < |\lambda| <
\frac{1}{\e_1}$ and $E \in \R$ (note that $L(E,\lambda f_i) \ge 1
$ if $|E| \ge  \| \lambda f_i \| + 4$). Our constructions start
with $F_1$ and we will divide them into
several steps. \\[1mm]
\noindent \textit{Construction 1.} First, we will apply Lemma~\ref{lem:a1}
to $F_1$ in order to enlarge the range of $\lambda$'s. Let $\e_2 =
\frac{min\{\e_1,\delta_1\}}{10}$. Then, there exists a finite family of
$\tilde{p}_1$-periodic potentials $\tilde{F}_1 = \{ \tilde{f}_1,\tilde{f}_2, \dots ,
\tilde{f}_{\tilde{m}_1} \} \subset B_{\e_0}(f)$ such that
$$
L(E,\lambda \tilde{F}_1) > \tilde{\delta}_1
$$
for some $0 < \tilde{\delta}_1 <1$ whenever $\forall \e_2 < |\lambda| < \frac{1}{\e_2}$ and $E \in \R$,  and
\begin{equation} \label{e1}
\left| L(E,\lambda \tilde{F}_1) - L(E,\lambda F_1) \right| < \frac{\e_2}{2}
\end{equation}
whenever $|E| < \frac{1}{\e_2}$ and $|\lambda| < \frac{1}{\e_2}$.

Explicitly, the construction of $\tilde{F}_1$  follows from the proof of \cite[Claim 3.1]{a}. For very large $\tilde{p}_1
> p_1$ (it must obey $p_1 | \tilde{p}_1$), choose $N_1(\tilde{p}_1)$ such that if $|E| < \frac{1}{\e_2}$, $|\lambda| \leq
\frac{1}{\e_2}$, $f_i \in F_1$ and a $\tilde{p}_1$-periodic potential $\tilde{f}$ which is $\frac{2 p_1 + 1}{N_1(\tilde{p}_1)}$
close to $f_i$ then $|L(E,\lambda \tilde{f}) - L(E,\lambda f_i)| < \frac{\e_2}{2}$, since the Lyapunov exponent is continuous for
periodic potentials (see Proposition \ref{indely}).

For $1 \leq j \leq 2 p_1 +1$, we define $\tilde{p}_1$-periodic
potentials $\tilde{f}^{(i,j)}$ by $\tilde{f}^{(i,j)}(n) = f_i(n),
0 \leq n \leq \tilde{p}_1-2$ and $\tilde{f}^{(i,j)}(\tilde{p}_1 -
1) = f_i{(\tilde{p}_1 - 1)} + \frac{j}{N_1(\tilde{p}_1)}$. By
\cite[Claim 3.4]{a}, there exists $j_0$ such that the spectrum of
$\tilde{f}^{(i,j_0)}$ has exactly $\tilde{p}_1$ components, that
is, all gaps of its spectrum are open. For convenience, we write
$\tilde{f}^{(i)} = \tilde{f}^{(i,j_0)}$. So there exists $h =
h(F_1,\tilde{p}_1,\e_2) > 0$ such that for any $f_i \in F_1$ and
$\e_2 \leq |\lambda| \leq \frac{1}{\e_2}$, $\Sigma(\lambda
\tilde{f}^{(i)})$ has $\tilde{p}_1$ components and the Lebesgue
measure of the smallest gap is at least $h$. Choose an integer
$N_2(\tilde{p}_1)$ with $N_2(\tilde{p}_1) > \frac{4 \pi }{\e_2 h
\tilde{p}_1}$.

For $0 \leq l \leq N_2(\tilde{p}_1)$, let $\tilde{f}^{(i,l)} =
\tilde{f}^{(i)} + \frac{4 \pi l }{\e_2 \tilde{p}_1
N_2(\tilde{p}_1)}$. Then $\tilde{F}_1$ is just the family obtained
by collecting the $\tilde{f}^{(i,l)}$ for different $f_i \in F_1$
and $0 \leq l \leq N_2(\tilde{p}_1).$ Order $\tilde{F}_1$ as
$\tilde{F}_1 = \{ \tilde{f}_1,\tilde{f}_2, \dots ,
\tilde{f}_{\tilde{m}_1} \} $ such that $\tilde{f}_1 =
\tilde{f}^{(1,0)}$ and $\tilde{f}_{\tilde{m}_1} =
\tilde{f}^{(1,1)}$. We can also assume that $N_2(\tilde{p}_1)$ was
chosen large enough, so that we have $\| \tilde{f}_{\tilde{m}_1} -
\tilde{f}_1 \| = \frac{4 \pi }{\e_2 \tilde{p}_1 N_2(\tilde{p}_1)}
< 1/3$ (this will be used to conclude that our limit-periodic
potentials are Gordon
potentials).\\[1mm]
\noindent \textit{ Construction 2}. Applying Lemma \ref{lem:a2} to $\tilde{F}_1$, there exists a finite family of $p_2$-periodic
potentials $ F_2 = \{ f^{\vec{t}_1}_2,f^{\vec{t}_2}_2, \dots, f^{\vec{t}_{m_2}}_2 \}$ such that
$$
F_2 \subset B_{p_2^{-2}} \subset B_{\e_2}  \subset B_{\e_0}(f)
$$
and
$$
L(E,\lambda F_2) > \delta_2
$$
for some $0 < \delta_2 < 1$ whenever $\e_2 < |\lambda| < \frac{1}{\e_2}$ and $E \in \R$, and
\begin{equation} \label{e2}
\left |L(E,\lambda F_2) - L(E,\lambda \tilde{F}_1)\right | < \frac{\e_2}{2}
\end{equation}
whenever $|E|$, $|\lambda| < \frac{1}{\e_2}$. From \eqref{e1} and \eqref{e2}, we have
$$
\left |L(E,\lambda F_2) - L(E,\lambda F_1)\right | < \e_2
$$
for $|E|$, $|\lambda| < \frac{1}{\e_2}$.

Explicitly, we construct $F_2$ as follows (cf.~the proof of
Lemma~\ref{lem:a2}). Let $p_2$ large and $p_2 = \tilde{m}_1
\tilde{p}_1 r_2 + d,\ 0 \leq d \leq \tilde{m}_1 \tilde{p}_1 - 1$.
Let $I_j = [j\tilde{p}_1,(j+1)\tilde{p}_1 -1] \subset \Z$ and let
$0 = j_0 < j_1 < \dots < j_{\tilde{m}_1-1} < j_{\tilde{m}_1} =
\frac{p_2}{\tilde{p}_1}$ be a sequence such that $j_{i+1} - j_i =
r_2 + 1$ when $0 \leq i < d/\tilde{p}_1$ and $j_{i+1} - j_{i} = r_2$ when $d/\tilde{p}_1
\leq i \leq \tilde{m}_1 \tilde{p}_1 - 1$. Define a $p_2$-periodic
potential $f_2(l)$ for $0 \leq l \leq p_2 - 1$ as follows. Let $j$
be such that $l \in I_j$ and let $i$ be such that $j_{i-1} \leq j
< j_i$ and let $f_2(l) = \tilde{f}_i(l).$ For any sequence
$\vec{t} = (t_1,t_2,\dots,t_{\tilde{m}_1})$ with $t_i \in
\{0,1,\dots,r_2 - 1\}$, let $f^{\vec{t}}_2 $ be a $p_2$-periodic
potential defined as follows. Let $0 \leq l \leq p_2 - 1$, and let
$j$ be such that $l \in I_j$. If $j = j_i - 1$ for some $1 \leq i
< \tilde{m}_1$, we let $f^{\vec{t}}_2(l) = f_2(l) + r^{-4}_2 t_i$,
and $j = j_{\tilde{m}_1} - 2$ then $f^{\vec{t}}_2(l) = f_2(l) +
r^{-4}_2 t_{\tilde{m}_1}$. Otherwise we let $f^{\vec{t}}_2(l) =
f_2(l)$. Let $p_2$ be sufficiently large so that $p^{-2}_2 < 1/3$.

Moreover, we can estimate the Lebesgue measure of the spectrum.
For any $E \in \R$ and $\e_2 < |\lambda| < \frac{1}{\e_2}$, we can
find $\tilde{f}_i \in \tilde{F}_1$ such that $L(E,\lambda
\tilde{f}_i) > \tilde{\delta}_1$ since $L(E,\lambda \tilde{F}_1)
> \tilde{ \delta}_1.$ If $r_2$ large enough, we have $\| A^{(E,
\lambda \tilde{f}_i)}_{(r_2 - 2) \tilde{p}_1} \| >
e^{\tilde{\delta}_1 (r_2 - 2) \tilde{p}_1}.$ Then we have
$$
\| A^{(E, \lambda f^{\vec{t}_k}_2)}_{(r_2 - 2) \tilde{p}_1}(f^{\vec{t}_k}_2(j_{i-1} \tilde{p}_1)) \| = \| A^{(E,
\lambda \tilde{f}_i)}_{(r_2 - 2) \tilde{p}_1} \| > e^{\tilde{\delta}_1 (r_2 - 2) \tilde{p}_1}.
$$
Since $E$ is arbitrary, we can apply Lemma \ref{lem:measure} to
conclude that the total Lebesgue measure of $\Sigma(\lambda
f^{\vec{t}_k}_2)$ is at most $4 \pi p_2 e^{- \tilde{\delta}_1 (r_2
- 2) \tilde{p}_1} < e^{-\tilde{p}_1 p^{1/2}_2}$ when $r_2$
sufficiently large. (Here
$f^{\vec{t}_k}_2$ can be any element from $F_2$.)\\[1mm]
\noindent \textit{ Construction 3}. Repeating the above procedures. Once we have constructed $F_{i-1} \subset B_{p^{-(i-1)}_{i-1}}
\subset B_{\e_{i-1}}$, by Lemma \ref{lem:a1}, we can get a finite family of $\tilde{p}_{i-1}$-periodic potentials $\tilde{F}_{i-1}
\subset B_{p^{-(i-1)}_{i-1}} $ satisfying the following. Let $\e_i = \frac{min\{\e_{i-1},\delta_{i-1}\}}{10}$, and we have
$$
L(E,\lambda \tilde{F}_{i-1}) > \tilde{\delta}_{i-1}
$$
for some $0 < \tilde{\delta}_{i-1} <1$ whenever $\forall \e_i < |\lambda| < \frac{1}{\e_i}$ and $E \in \R$,  and
$$
\left| L(E,\lambda \tilde{F}_{i-1}) - L(E,\lambda F_{i-1}) \right| < \frac{\e_i}{2}
$$
whenever $|E| < \frac{1}{\e_i}$ and $|\lambda| < \frac{1}{\e_i}$.

Next, as in \textit{Construction 2}, we will get a finite family
$F_i$ of $p_i$-periodic potentials which satisfies the following
(here our perturbation is $r^{-N i}_i t = r^{-2 i}_i t$, $t \in \{0, 1, 2, \dots, r_i - 1 \}$). \\[1mm]
(i). $L(E,\lambda F_i) > \delta_i$ for some $0 < \delta_i < 1$ and all $E \in \R$ and $\e_i < |\lambda| < \e_i^{-1}$.\\
(ii). $\left |L(E,\lambda F_i) - L(E,\lambda F_{i-1}) \right | < \e_i$, for $|E| < \frac{1}{\e_i}$ and
$|\lambda| < \frac{1}{\e_i}$.\\
(iii). $F_i \subset B_{p^{-i }_i} \subset B_{\e_i} \subset B_{\e_{i-1}}\subset B_{\e_0}(f), i > 2.$ (Note: $B_{\e_2}$ may not be
in $B_{\e_1}$.)\\ (iv). $\forall f^{\vec{t}_k}_i \in F_i$, $|\Sigma(\lambda f^{\vec{t}_k}_i)| \leq e^{-\tilde{p}_{i-1} p^{1/2}_i}$ when $\e_i <
|\lambda| < \e_i^{-1}$ (here $|\cdot|$ denotes the Lebesgue measure).\\
(v). $ \label{ineq:g}p^{- i}_i < \frac{1}{3} (i-1)^{-\tilde{p}_{i-1}}$ since we can let $p_i$ be sufficiently large.\\
(vi). $ \| f^{\vec{t}_1}_i - f^{\vec{t}_{m_i}}_i \| = \frac{4
\pi}{\e_i \tilde{p}_{i-1} N_2(\tilde{p}_{i-1})} < \frac{1}{3}
(i-1)^{-\tilde{p}_{i-1}}.$  Here $N_2(\tilde{p}_{i-1})$ appears as
in \textit{Construction 1}, and we can ensure that this inequality
holds since
$\tilde{p}_{i-1}$ is fixed while $N_2(\tilde{p}_{i-1})$ can be taken as large as needed.\\[1mm]
\indent Then we will get a limit-periodic potential $f_{\infty}
\in B_{\e_0}(f)$, whose Lyapunov exponent is a positive continuous
function of energy $E$ and the Lebesgue measure of the spectrum is
zero (Lemma \ref{lem:limitlya} implies that $L(E,\lambda
f^{\vec{t}}_i) \rightarrow L(E,\lambda f_\infty)$). Moreover, we
have the following two claims.
\begin{claim} \label{claim:gon}
$f_{\infty}$ is a Gordon potential.
\end{claim}

\begin{proof}
Let $q_i = \tilde{p}_i$. Obviously, $q_i \rightarrow \infty$ as $i
\rightarrow \infty.$  For $i \ge 1$, we have
\begin{align*}
\max_{1\leq n \leq q_i}\left|f_{\infty}(n)-f_{\infty}(n\pm q_i)
\right| &\leq |f_{\infty}(n) - f^{\vec{t}_1}_{i+1}(n) | +
|f_{\infty}(n\pm q_i) - f^{\vec{t}_1}_{i+1}(n\pm q_i) | \\& \quad +
|f^{\vec{t}_1}_{i+1}(n) - f^{\vec{t}_1}_{i+1}(n\pm q_i) |\\& \leq
p^{- (i+1)}_{i+1} + p^{- (i+1)}_{i+1} + \frac{4 \pi }{\e_{i+1}
\tilde{p}_{i} N_2(\tilde{p}_{i})}\\& \leq 2 \frac{1}{3}
(i)^{-\tilde{p}_{i}} + \frac{1}{3} (i)^{-\tilde{p}_{i}} \\& \leq
i^{-q_i}.
\end{align*}
So $f_{\infty}$ is a Gordon potential. (Here $f^{\vec{t}_1}_{i+1}$
is an element of $F_{i+1}$).
\end{proof}

\begin{claim} \label{claim:haus}
$\Sigma(\lambda f_{\infty})$ has zero Hausdorff dimension for every $\lambda \not= 0$.
\end{claim}

\begin{proof} Let $\lambda \neq 0$ and $0 < \alpha \leq 1$ be given. Without loss of generality,
assume $\lambda > 0$. Choose $i$ large enough so that $\e_i < \lambda < 1/\e_i$ and $1/i < \alpha$.
For every $f^{\vec{t}_k}_i \in F_i$,  $\| \lambda f_{\infty} - \lambda f^{\vec{t}_k}_i \| < \lambda p^{-i}_i$
implies\footnote{It is well known that for $V,W : \Z \to \R$ bounded, we have $\mathrm{dist} ( \sigma(\Delta+V) ,
\sigma(\Delta+W) ) \le \| V - W \|_\infty$, where $\mathrm{dist}(A,B)$ denotes the Hausdorff distance of
two compact subsets $A , B$ of $\R$.} $\mathrm{dist} (\Sigma(\lambda f_{\infty}),\Sigma(\lambda f^{\vec{t}_k}_i))
< \lambda p^{-i}_i.$ Since $\lambda f^{\vec{t}_k}_i$ is $p_i$-periodic, we have
$$
\Sigma(\lambda f^{\vec{t}_k}_i) = \bigcup^{p_i}_{z=1} \tilde{I}^{(\vec{t}_k,i)}_z,
$$
where $\tilde{I}^{(\vec{t}_k,i)}_z = [a_z,b_z]$ is a closed interval.

Let $I^{(\vec{t}_k,i)}_z = [a_z - \lambda p^{-i}_i,b_z + \lambda p^{-i}_i]$ and since $\mathrm{dist}
(\Sigma(\lambda f_{\infty}),\Sigma(\lambda f^{\vec{t}_k}_i)) \leq \lambda p^{-i}_i$, we have
$$
\Sigma(\lambda f_\infty) \subset \bigcup^{p_i}_{z=1} I^{(\vec{t}_k,i)}_z.
$$
Moreover, $b_z - a_z \leq e^{-\tilde{p}_{i-1} p^{1/2}_i}$ since $|\Sigma(\lambda f^{\vec{t}_k}_i)|
\leq e^{-\tilde{p}_{i-1} p^{1/2}_i}$. Then we have
\begin{align*}
h^{\alpha}(\Sigma(\lambda f_{\infty})) & \leq \lim_{i\rightarrow \infty}\sum^{p_i}_z (e^{-\tilde{p}_{i-1}
p^{1/2}_i} + 2 \lambda p^{-i}_i)^{\alpha}\\
& = \lim_{i\rightarrow \infty}p_i (e^{-\tilde{p}_{i-1} p^{1/2}_i} + 2 \lambda p^{-i}_i)^{\alpha} \\
& = \lim_{i\rightarrow \infty}(p^{1/\alpha}_i e^{-\tilde{p}_{i-1} p^{1/2}_i} + 2 \lambda p^{-i + 1/\alpha}_i)^{\alpha}.
\end{align*}
Since $1/i < \alpha$, we have $-i + 1/\alpha < 0$, and it follows that
$$
\lim_{i\rightarrow\infty} (p^{1/\alpha}_i e^{-p_{i-1} p^{1/2}_i} + 2 \lambda p^{-i  +1/\alpha}_i)^{\alpha} = 0.
$$
So we have $ h^{\alpha}(\Sigma(\lambda f_{\infty})) = 0 $ (note:
when $i \rightarrow \infty$, $\lambda$ belongs to
$(\e_i,\frac{1}{\e_i})$ for all $i$ large enough since this
interval is expanding). So the Hausdorff dimension of the spectrum
is less than $\alpha$. Since $\alpha$ was arbitrary, the Hausdorff
dimension must be zero.
\end{proof}

This implies all the assertions in Theorem~\ref{t.main} except for
the absence of eigenvalues for every $\omega$. Given the Gordon
Lemma (see Lemma~\ref{lem:gordon} above), this last statement will
follow once Theorem~\ref{t.gordon} is established.
\end{proof}

\begin{remark}
Since $\delta_i \leq \delta_{i-1}/10, i \ge 1$, it is true that when $\e_i < |\lambda| < \frac{1}{\e_i}$, $L(E,\lambda f_\infty) \ge
\frac{8}{9} \delta_i$ for any $E \in \R$. This gives information about the range of the Lyapunov exponent on certain intervals. Clearly, $\delta_i
\rightarrow 0$ when $i \rightarrow \infty$ since the Lyapunov exponent will go to zero when $\lambda$ goes to zero.
\end{remark}

\begin{proof}[Proof of Theorem~\ref{t.gordon}]
Let $\omega = e$ first. Relative to any minimal translation $\tilde{T}$, the selected $f$ in the proof of
Theorem~\ref{t.main} is still $n_0$-periodic by Proposition \ref{inde}, so we can start with the same ball
$B_{\e_0}(f)$ and choose the same periodic potentials in $B_{\e_0}(f)$. Then we get the same $f_\infty$. For
the finite family $F_i$ from \textit{Construction 3}, though the Lyapunov exponent may change, the following
properties hold (note that $\| f^{\vec{t}_1}_i \|$ does not change).\\[1mm]
(i). $F_i \subset B_{p^{-i }_i} \subset B_{\e_i} \subset B_{\e_0}(f).$\\[1mm]
(ii). $p^{- i}_i < \frac{1}{3} (i-1)^{-\tilde{p}_{i-1}}$.\\[1mm]
(iii). $\| f^{\vec{t}_1}_i - f^{\vec{t}_{m_i}}_i \| = \frac{4 \pi}{\e_i \tilde{p}_{i-1} N_2(\tilde{p}_{i-1})} < \frac{1}{3}
(i-1)^{-\tilde{p}_{i-1}}$.\\[1mm]
Then Claim~\ref{claim:gon} holds true, and so $f(\tilde{T}^n(e))$ is a Gordon potential. For arbitrary $\tilde{\omega}$, if we repeat the same
procedures, (i)---(iii) above still hold as stated (since none of them are related to $\omega$), and Theorem~\ref{t.gordon} follows.
\end{proof}

\end{document}